\documentclass[10pt, a4paper, reqno]{amsart}

\usepackage{amscd}
\usepackage{dsfont}
\usepackage{booktabs}
\usepackage{changepage}

\usepackage[indentafter]{titlesec}
\titleformat{name=\section}{}{\thetitle.}{0.8em}{\centering\scshape}
\titleformat{name=\subsection}{}{\thetitle.}{0.8em}{\centering\scshape}

\usepackage{wasysym} 
\usepackage[utf8]{inputenc}

\usepackage{pgfplots}
\pgfplotsset{width=7cm,compat=1.8}
\usepackage{pgfplotstable}

\usepackage{verbatim}

\usepackage{amsmath}
\usepackage{amstext}
\usepackage{amsthm}
\usepackage{amssymb}
\usepackage{amsopn}
\usepackage{amssymb}
\usepackage{amsxtra}
\usepackage{amsfonts}
\usepackage{fancyhdr, xcolor}
\usepackage{paralist, epsfig}
\usepackage[mathcal]{euscript}
\usepackage[english]{babel}

\makeatletter
\g@addto@macro\normalsize{%
  \setlength\abovedisplayskip{10pt}
  \setlength\belowdisplayskip{10pt}
  \setlength\abovedisplayshortskip{5pt}
  \setlength\belowdisplayshortskip{8pt}
}

\usepackage{pgfplots}
\pgfplotsset{compat=newest}

\usepgfplotslibrary{patchplots}

\theoremstyle{definition}
\newtheorem{h}{Problem}

\usepackage{setspace}
\usepackage{color}

\definecolor{grey}{gray}{.3}

\usepackage{hyperref}
\hypersetup{
    linkbordercolor = {red},
}
\newcommand\blfootnote[1]{%
  \begingroup
  \renewcommand\thefootnote{}\footnote{#1}%
  \addtocounter{footnote}{-1}%
  \endgroup
}

\begin{document}

$ $

\vspace{-40pt}

\title{Engaging students in conjecturing through homework\\in Real Analysis and Differential Equations}

\author{Thomas Pawlaschyk\hspace{0.5pt}\MakeLowercase{$^{\text{a}}$} and Sven-Ake Wegner\hspace{0.5pt}\MakeLowercase{$^{\text{b}}$}}

\hspace{-1000pt}\blfootnote{\hspace{5.5pt}2010 \emph{Mathematics Subject Classification}: Primary 97C70; Secondary 97C30, 97D40, 97B40\vspace{1.6pt}}

\hspace{-1000pt}\blfootnote{\hspace{5.5pt}\emph{Key words}: Real Analysis, Differential Equations, Conjecture, Refutation, Proof, Counterexample, Task Design.\vspace{1.6pt}}


\hspace{-1000pt}\blfootnote{\hspace{0pt}$^{a}$\,University of Wuppertal, School of Mathematics and Natural Sciences, Gau\ss{}stra\ss{}e 20, 42119 Wuppertal, Germany,\linebreak\phantom{x}\hspace{1.2pt}phone: +49\,(0)\,202\:439\:25\:41, e-mail: pawlasch@uni-wuppertal.de.\vspace{1.6pt}}

\hspace{-1000pt}\blfootnote{\hspace{0pt}$^{b}$\,Corresponding author: Teesside University, School of Science, Engineering and Design, Middlesbrough, TS1\;3BX,\linebreak\phantom{x}\hspace{1.2pt}United Kingdom; currently: Institut des Hautes \'Etudes Scientifique, Universit\'e Paris-Saclay, 91440 Bures-sur-\linebreak\phantom{x}\hspace{1.2pt}Yvette, France, phone: +33\,(0)\,160\:92\:66\:61, e-mail: wegner@ihes.fr.\vspace{3pt}}

\begin{abstract} In this note we report on an implementation of discovery-oriented problems in courses on Real Analysis and Differential Equations. We explain a type of task-design that gives students the opportunity to conjecture, refute and prove. What is new is that the complexity in our problems is limited and thus the tasks can also be used in homework assignments. In addition to several concrete examples we also discuss feedback and assessment outcomes of our students.
\end{abstract}

\maketitle

\vspace{-15pt}


\section{Introduction}\label{SEC:1}

There is a wide consensus in the mathematics education community that activities should be implemented in the practice of teaching that engage students in cognitively demanding tasks which go beyond carrying out computations according to memorized procedures in an algorithmic fashion \cite{HK18, KoJo2019, Li}. These activities should encourage students to develop mathematical creativity \cite{Li}, facilitate conceptual knowledge \cite{HK18, K2008} and improve the ability to transform a given collection of information into a well-posed mathematical question such that then mathematical tools can be employed to gain a solution \cite{W1980}. Formulated more generically, these activities should be closer to what mathematicians do and devolve `mathematical responsibility' to students \cite{D14, F2016, K2008, W1980}.

\smallskip

Up to now there are two, basically different but not conflicting, major approaches to implement this idea. The first is to engage students either in an actual research project \cite{PrimusI, PrimusII} or to emulate such research projects \cite{A77, FM17, W1980}. Here, students work outside of regular lessons, participation can be voluntarily and assessment might be dropped\footnote{\,In some education systems, e.g., in the UK or in Germany, BSc programs include a final year project, in which\linebreak\phantom{x}\hspace{1.2pt}students engage in (emulated) research activities. In this context participation is of course mandatory and assessment\linebreak\phantom{x}\hspace{0.8pt}takes place.}. The second approach is to incorporate activities as mentioned in the first paragraph into regular lessons. This approach thus comes with the following constraints. Topics have to be selected according to a prescribed curriculum, activities have to fit into a given timetable, all students have to be included in the activity, and assessment has to be carried out in a reasonable way taking into account the institution's conventions.

\smallskip

In this paper we report on an implementation of the second approach in courses on Real Analysis and Differential Equations taught at a German university in 2016-17. In both cases the structure of the course involved two classical 90-minutes blackboard lectures per week, one 90-minutes recitation, weekly homework and a final exam, the latter written in RA and oral in DE. For this reason we sought to implement our activities in form of special homework problems and tasks to be done during the recitations.

\smallskip

In Section \ref{SEC:2} we give a short overview of similar activities that have been surveyed in the recent literature. In Section \ref{SEC:3} we explain our setup and explain a task design that we developed to foster our students' engagement in particular with conjectures. This special task design seems to our knowledge not to appear in the literature so far. In Section \ref{SEC:4} we present feedback of our students.


\bigskip
\section{Discovery-oriented activities in the classroom}\label{SEC:2}
\smallskip

When it comes to higher mathematics education, many teachers follow a `definition-theorem-proof' scheme. That is, abstract definitions are followed by formal theorems that relate the mathematical objects or properties that have been defined. Then proofs for the theorems are given. Many teachers enrich this approach via remarks, that explain why definitions are made as they are made, or by outlining which alterations of theorems are possible and which are not. Both often includes examples, non-examples and counterexamples. These illustrate for instance that a theorem is sharp or that a definiton rules out a certain pathology. Consequently, traditional tasks like

\vspace{6pt}

\begin{compactitem}

\item[$\bullet$] `Compute [\,\dots\,].'

\vspace{5pt}

\item[$\bullet$] `Solve [\,\dots\,].'

\vspace{5pt}

\item[$\bullet$] `$\frac{\operatorname{d}}{\operatorname{d}\hspace{-1pt}x}[\,\dots\,]=\,$?'

\end{compactitem}

\vspace{7pt}

on which students work during recitations, as homework, or in exams are (partly) replaced by tasks like

\vspace{0pt}

\begin{compactitem}

\item[$\bullet$] `Prove [\,\dots\,].'\vspace{5pt}

\item[$\bullet$] `Find an example [\,\dots\,].'\vspace{5pt}

\item[$\bullet$] `Prove or disprove [\,\dots\,].'\vspace{7pt}

\end{compactitem}

Based on this, but pushing the boundary further to achieve what we mentioned at the beginning of Section \ref{SEC:1}, tasks like

\vspace{6pt}

\begin{compactitem}

\item[$\bullet$] `Find a formula for [\,\dots\,].' \cite[p.\ 941]{Li} and \cite[p.\ 163]{W1980}\,\footnote{\,Both tasks come with a picture of a graph $\Gamma_r$ depending on a parameter $r$. Lithner's \cite{Li} task is to find the number\linebreak\phantom{x}\hspace{1.2pt}of edges depending on $r$.  Watson's \cite{W1980} is more complicated.\vspace{2pt}}

\vspace{5pt}

\item[$\bullet$] `Complete in a way that you believe makes true statements: [\,\dots\,].' \cite[p.\ 742]{F2016}

\vspace{5pt}

\item[$\bullet$] `Generalize [\,\dots\,].' \cite[\S\,I.2, Exercise 16(d) on p.~A\,I.121]{B2}\,\footnote{\,Bourbaki \cite{B2} requires first to prove some proposition and then to employ the latter to generalize a theorem that was\linebreak\phantom{x}\hspace{1.2pt}given in a previous chapter.\vspace{2pt}}

\vspace{5pt}

\item[$\bullet$] `Define [\,\dots\,].' \cite[p.\ 94]{D14} and \cite[p.\ 738]{F2016}\,\footnote{\,Dawkings' \cite{D14} task is to `re-invent' the notion of a group based on examples (real numbers, integers, invertible\linebreak\phantom{x}\hspace{1.2pt}matrices) in which the equations $a*x=b$ with $*\in\{\cdot_\mathbb{R},+_\mathbb{Z},\cdot_{\text{GL}_n}\}$ can be solved. Fukawa-Connelly \cite{F2016} asked\linebreak\phantom{x}\hspace{1.2pt}students to find a definition for $\lim_{x\rightarrow p}f(x)=\infty$ after the definition of $\lim_{x\rightarrow p}f(x)=L$ for $L\in\mathbb{R}$ was discussed\linebreak\phantom{x}\hspace{1pt}in a previous lesson.}

\vspace{5pt}

\item[$\bullet$] `Prove or refute [\,\dots\,].' or with more words `Prove or demonstrate a counterexample for the following propositions (if you can `save' an incorrect proposition by a slight modification---say barring a certain example, do so): [\,\dots\,].' \cite[p.\ 444]{K2008} and \cite[p.\ 742]{F2016}

\vspace{5pt}

\end{compactitem}

\vspace{2pt}

have been suggested in the literature. And even more discovery-oriented activities can be found. In \cite{K2008} the `If not, what yes?' approach (hilarious, a must-read for every math professor!) starts with a plausible but wrong conjecture that is presented by the teacher and leads to a chain of counterexamples and refutations until in the end a correct statement is proved. In this way the atmosphere of a `faithful' research activity \cite[p.\ 452]{K2008} is created in the classroom. In \cite{LL} a detailed description of how students can `re-invent' the concept of quotient groups is given, and in \cite{KoJo2019, KKTS2014} geometric theorems have to be discovered and proved.

\smallskip

A central aim in all aforementioned articles is to put the learners in the center of activity and to devolve `mathematical responsibility' to them. In order to achieve this, many classroom setups require however that the teacher at least assumes the role of a moderator during student-led discussions. Sometimes it is even necessary to intervene when the discussion is stuck or goes off-topic. Thus, for a given course, e.g., with a classical structure of weekly lectures and recitations, graded homework, and written exams, not all of the concepts above are applicable. In \cite[p.\ 743]{F2016} the author for instance states that `there were no instances in homework or exams where the professor asked the students to engage in conjecturing activities.'

\smallskip

The two courses, that we are going to expand on in Section \ref{SEC:3} and that led to the current paper, were prepared with the aim to implement discovery-oriented activities into homework and also into exams. For this reason we were looking for a way to keep the discovery element in the problem, but at the same time to limit the `degrees of freedom' such that the problem becomes suitable as a homework or exam task.

\smallskip

In the sequel we adopt Koichi's notation \cite{K2008} and refer to all problems of the type characterized in this section and in Section \ref{SEC:1} as \emph{discovery-oriented problems}.


\section{A task design to foster conjecturing}\label{SEC:3}
\smallskip

\subsection{Ordinary Differential Equations}
\smallskip

The starting point for this article was a course given in the winter term 2016-17 at the University of Wuppertal (Germany) about Ordinary Differential Equations. This course gave an introduction to basic existence, uniqueness and extension theory for solutions, stability theory and solution methods. Our students mostly were mathematics majors in the third year of their Bachelor studies. The assessment consisted of weekly homework sheets and an oral exam at the end of the course. The final grade for the course was the grade in the oral exam, but raised by one step (e.g., A- instead of B+) if the student scored between 50\% and 70\% on the homework sheets, and by two steps (e.g., A instead of B+)  if the student scored above 70\% on the homework sheets. When designing the homework tasks, our aim consequently was to provide a fair chance to score above 50\% resp.~70\% and at the same time to prepare for the oral exam. In the latter, questions focus traditionally more on theoretical understanding than on computations.

\smallskip

In particular to engage students in conjecturing and refuting, we used the following design to create tasks that are on the one hand `open' in the sense of Section \ref{SEC:2} but on the other hand can be used as (assessed) homework problems in the sense that they are `doable' for the students under the given constraints in time, background, ability etc. Our approach is best explained with the following prototype example that we employed in order to clarify the assumptions of the Picard-Lindel\"of theorem on existence and uniqueness of solutions of the initial value problem $y'=f(x,y)$, $y(x_0)=y_0$.

\vspace{-5pt}

\begin{adjustwidth}{1cm}{1cm}
\begin{h}\label{P-1} Let $f\colon Q:=[a,b]\times[c,d]\rightarrow\mathbb{R}$, $(x,y)\mapsto f(x,y)$ be a continuous map. Compare the following three conditions.\vspace{5pt}
\begin{compactitem}
 \item[(i)] $f$ is Lipschitz continuous.\vspace{5pt}
 \item[(ii)] $f$ satisfies a Lipschitz condition with respect to $y$.\vspace{5pt}
  \item[(iii)]$f(x,\cdot)\colon[c,d]\rightarrow\mathbb{R}$ is Lipschitz continuous for every $x\in[a,b]$.\hfill$\diamondsuit$
\end{compactitem}
\end{h}
\end{adjustwidth}

\medskip\smallskip

We recall that $f$ as above is Lipschitz continuous, if
\begin{equation}\label{LIP}
\exists\:L\geqslant0\:\forall\:(x_1,y_1),(x_2,y_2)\in Q\colon |f(x_1,y_1)-f(x_1,y_2)|\leqslant L\|(x_1,y_1)-(x_2,y_2)\|
\end{equation}
holds, that $f$ satisfies a Lipschitz condition with respect to $y$ if
\begin{equation}\label{LIP-1}
\exists\:L\geqslant0\:\forall\:(x,y_1),(x,y_2)\in Q\colon |f(x,y_1)-f(x,y_2)|\leqslant L|y_1-y_2|
\end{equation}
holds, and that $f(x,\cdot)\colon[c,d]\rightarrow\mathbb{R}$ is Lipschitz continuous for every $x\in[a,b]$, if
\begin{equation}\label{LIP-2}
\forall\:x\in[a,b]\:\exists\:L\geqslant0\:\forall\:y_1,\,y_2\in [c,d]\colon |f(x,y_1)-f(x,y_2)|\leqslant L|y_1-y_2|
\end{equation}
holds. In \eqref{LIP} the symbol $\|\cdot\|$ denotes an arbitrary norm on $\mathbb{R}^2$. Looking at the quantifiers it is clear that (i)$\Rightarrow$(ii)$\Rightarrow$(iii) can be seen by putting $x=x_1=x_2$ and $L(x)=L$, respectively. In order to see that (ii)$\Rightarrow$(i) is not valid one can adapt the standard counterexample that shows that on a compact interval continuity and Lipschitz continuity are not the same. More precisely, one defines $f\colon[0,1]\times[0,1]\rightarrow\mathbb{R}$, $f(x,y)=\sqrt{x}$. Constructing a counterexample for (iii)$\Rightarrow$(ii) is probably the hardest part of Problem \ref{P-1}. One approach is to look for a function that has at all points $(x,0)$ in $y$-direction a finite gradient but the gradients tend to infinity when $x$ approaches zero, see Figure 1 for a picture. An explicit counterexample is for instance $f\colon[0,1]\times[0,1]\rightarrow\mathbb{R}$, $f(x,y)=(x+y)^{1/2}x^{1/3}$.

\bigskip

\begin{center}

\begin{tikzpicture}[scale=1.25]
    \begin{axis}[
    xtick = {1, 0.75, 0.5, 0.25, 0},
    ytick = {0, 0.25, 0.5, 0.75, 1},
    ztick = {0, 0.3, 0.6, 0.9, 1.2},
    zmin=-0.01,zmax=1.5,
    yticklabels=\empty,
    xticklabels=\empty,
    zticklabels=\empty,
	grid,
	colormap={summap}{
        color=(red!80); color=(red!80); color=(red!80)
    }
]
    \addplot3[color=red, patch,patch refines=4,
		shader=faceted,
		patch type=biquadratic] 
    table[z expr=(1-x+y)^(1/2)*(y)^(1/100)]
    {
        x  y
        0 0
        1  0
        1  1
        0 1
        0.5  0
        1  0.5
        0.5  1
        0 0.5
        0.5  0.5
    };
    \end{axis}
    
\node at (0, 0.4)   {$1$};

\node at (1.7, 0.075)   {$y$};

\node at (3.725, -0.23)   {$0$};

\node at (4.65, 0.4)   {$x$};

\node at (5.425, 1.05)   {$1$};

\end{tikzpicture}

\medskip

{

\small

\textbf{Figure 1.}\label{F-1}~Plot of the function $f\colon[0,1]^2\rightarrow\mathbb{R}$, $f(x,y)=(x+y)^{1/2}x^{1/3}$.

}

\end{center}

\smallskip

During the discussions with our students we observed immediately, that the formulation `Compare the conditions [\,\dots\,]' in Problem \ref{P-1} appeared to be very unusual to them and that it was hard for them to understand what was required to be done in Problem \ref{P-1}. Indeed, the students were used to formulations like `Show [\,\dots\,]', `Prove [\,\dots\,]' or `Find [\,\dots\,]' followed by detailed instructions. We give a classical version of Problem \ref{P-1} below.

\begin{adjustwidth}{1cm}{1cm}
\begin{h}\label{P-2} Let $f\colon Q:=[a,b]\times[c,d]\rightarrow\mathbb{R}$, $(x,y)\mapsto f(x,y)$ be a continuous map. 
\vspace{3pt}
\begin{compactitem}
 \item[(i)] Assume that $f$ is Lipschitz continuous. Prove that $f$ satisfies a Lipschitz condition with respect to $y$.
 
 \vspace{3pt}

\item[(ii)] Assume that $f$ satisfies a Lipschitz condition with respect to $y$. Prove that the map $f(x,\cdot)\colon[c,d]\rightarrow\mathbb{R}$ is Lipschitz continuous for every $x\in[a,b]$.

\vspace{3pt}

\item[(ii)] Let $f\colon[0,1]\times[0,1]\rightarrow\mathbb{R}$, $f(x)=(x+y)^{1/2}x^{1/3}$. Show that $f(x,\cdot)$ is Lipschitz continuous for every $x\in[0,1]$. Show that $f$ does not satisfy a Lipschitz condition with respect to $y$.

\vspace{3pt}

\item[(iv)] Find a function $f\colon[0,1]\times[0,1]\rightarrow\mathbb{R}$, $(x,y)\mapsto f(x,y)$ that satisfies a Lipschitz condition with respect to $y$ but which is not Lipschitz continuous.\hfill$\diamondsuit$

\end{compactitem}
\end{h}

\end{adjustwidth}

\bigskip

It is clear that Problem \ref{P-2} is easier in that one of the two counterexamples is given. In addition it is considerably harder to prove the correct implications and find counterexamples for the wrong implications if one does not know a priori which ones are correct and which ones are wrong. However, since there are only three conditions given, the number of conjectures to start with is limited and the setup allows only for one refutation per conjecture as the next refutation would lead back to the original conjecture. This design thus gives also to those who have a priori no intuition, the chance to start by just trying to prove one implication selected at random or to check all three conditions in one random example. And indeed, this might be exactly what a mathematician would do in a scientific context if she encounters a new and unfamiliar problem for the first time.

\smallskip

Of course there are versions of Problem \ref{P-1} that devolve much more responsibility to the students, see for instance the following variant.

\begin{adjustwidth}{1cm}{1cm}
\begin{h}\label{P-3} Discuss the notion of a Lipschitz condition for $f\colon [a,b]\times[c,d]\times\mathbb{R}\rightarrow\mathbb{R}$, $(x,y)\mapsto f(x,y)$ with respect to $y$.\hfill$\diamondsuit$
\end{h}
\end{adjustwidth}

\bigskip

This exercise would include firstly the task to find reasonable other conditions that one can compare with the condition to be discussed. A problem like this goes more in the direction of undergraduate research, see Section \ref{SEC:1}, or would need an effective moderation by the teacher, see Section \ref{SEC:2}. Consequently, we did not use Problem \ref{P-3} in the course. However, as a follow-up, we discussed in a Bachelor thesis, see \cite{PW}, that changing the quantifiers in the Lipschitz condition, i.e., replacing \eqref{LIP-1} with \eqref{LIP-2}, leads not to a sufficient condition for the Picard-Lindel\"of theorem to hold.

\smallskip

\subsection{Real Analysis II}

This subsection is devoted to another example of a problem that we used in a course taught at the University of Wuppertal (Germany) in the summer term 2017. The course Real Analysis II is the middle part of the sequence of courses Real Analysis I-III which are all mandatory for mathematics majors and which cover continuity, differentiation and integration of functions first in one and then in several real variables. In this case the assessment was a written exam and the students were required to obtain more than 50\% overall on the homework sheets to be admitted to this exam. There was no option to improve the exam grade through the exercises this time. Our prototype task is again to determine the relation between central concepts introduced in the lecture\footnote{\,Indeed, questions of the form `Which relationship can you find between [\,\dots\,]?' have been suggested in \cite[p.\ 322]{FOP2001}\linebreak\phantom{x}\hspace{1pt}but no concrete example was given and the article suggests that the term `relationship' refers rather to geometric\linebreak\phantom{x}\hspace{-0.4pt}`configurations' than to implications between (abstract) properties.}.

\begin{adjustwidth}{1cm}{1cm}
\begin{h}\label{P-4} Let $U\subseteq\mathbb{R}^n$ be open and $x_0 \in U$. For $f\colon U \to \mathbb{R}^m$ we consider the statements (i)--(vi) below. What is the relation of these statements? Justify your answer.\vspace{3pt}
\begin{compactitem}
\item[(i)] $f$ is continuously partially differentiable at $x_0$.\vspace{3pt}
\item[(ii)] $f$ is totally differentiable at $x_0$.\vspace{3pt}
\item[(iii)] $f$ is continuous and partially differentiable at  $x_0$.\vspace{3pt}
\item[(iv)] $f$ is partially differentiable at $x_0$.\vspace{3pt}
\item[(v)] $f$ is continuous at $x_0$.\vspace{3pt}
\item[(vi)] There is a unique matrix $A \in \mathbb{R}^{m \times n}$ such that
$$
{\displaystyle\lim_{\substack{x\to x_0 \\ x\neq x_0}}}\frac{f(x)-f(x_0)-A(x-x_0)}{|x-x_0|}=0
$$
holds.\hfill$\diamondsuit$
\end{compactitem}
\end{h}
\end{adjustwidth}

\bigskip

For the solution the students firstly needed to observe that (i)$\Rightarrow$(ii)$\Rightarrow$(iii)$\Rightarrow$(iv) and (iii)$\Rightarrow$(v) had been shown already in the lecture or follow directly from the definition. Secondly, they had to provide counterexamples in order to show that none of these implications is an equivalence. Then the students had to prove that (ii) and (vi) are equivalent and finally they had to show that neither (iv)$\Rightarrow$(v) nor (v)$\Rightarrow$(iv) is true, again by giving suitable counterexamples. Figure 2 contains the full picture of all valid and non valid implications.

\begin{center}
\includegraphics[width=115pt]{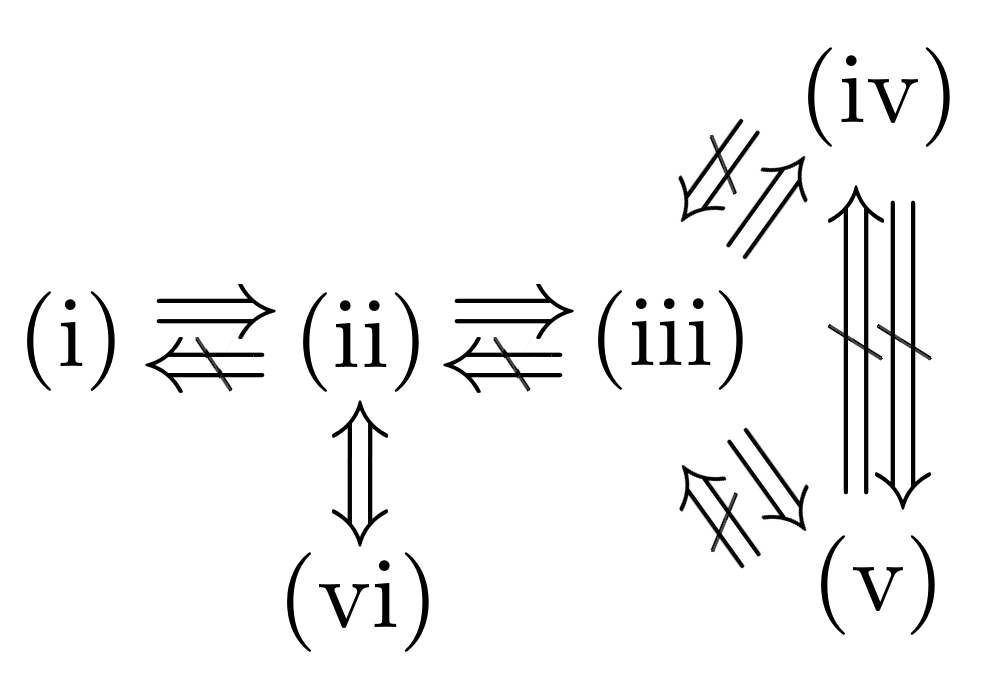}

{

\small

\textbf{Figure 2.}~Scheme of implications corresponding to Problem \ref{P-4}.

}

\end{center}

\medskip

Thirty students worked on this problem. The first part, citing corresponding results from the lecture, was indeed done correctly by 26. Only three students gave the counterexamples that we mentioned above as the second step, while the rest of the students did not discuss the converse implications at all. Only six students properly showed that (ii) and (vi) are equivalent. The other students only mentioned that (ii)$\Rightarrow$(vi) follows directly from the definition but did not give further arguments. In particular, they did not show that (ii) determines the matrix $A$ uniquely. The last point, the relation of (v) and (iv) was not discussed by any student. We mention that all students ignored that most of the statements from the lecture were only stated for functions $f\colon U\to\mathbb{R}$ but that here vector-valued functions appear which requires to use the arguments of the lecture for each coordinate function.

\medskip

The grading result of Problem \ref{P-4} is given in the following table.

\medskip

\begin{center}

{

\small

\textbf{Table 1.}\label{T-1}~Grading results for Problem \ref{P-4}.

}
\medskip

\begin{tabular}{lccc}
    \toprule
         & N\hspace{-3pt}$\phantom{\displaystyle X^{X^X}}$ & \hspace{10pt}Mean\hspace{10pt} & \hspace{22pt}SD\vspace{2pt}\\ \midrule\vspace{3pt}
\hspace{-5pt}\begin{minipage}{175pt}Grading results for Problem \ref{P-4}\phantom{$X^{X^{X}}$}\\\hspace{3pt}(min=0 points, max=10 points)\end{minipage}\hspace{5pt}\hspace{0pt} & \hspace{-21pt}30 & 2.5 &  \hspace{22pt}2.0 \\\bottomrule
  \end{tabular}
\end{center}

\bigskip

Although the format of the problem is exactly the same as in Problem \ref{P-1}, the complexity is much higher since the number of conditions has doubled. The grading results indicate that this level of complexity was maybe too high. There are however other factors, e.g., that Problem \ref{P-4} was only one out of five tasks for one week and that it was---compared to the first course that we mentioned in this section---not possible to raise the exam grade by a high score on the homework sheets. We conclude with the following solution handed in by one of the students.

\medskip

\begin{center}

\includegraphics[width=240pt]{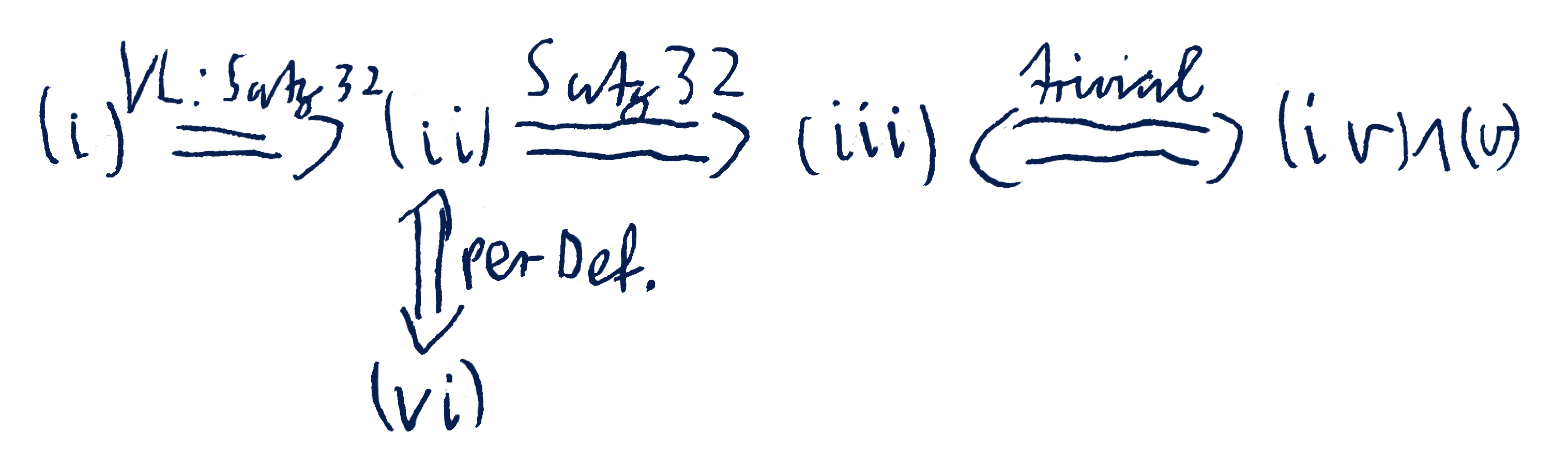}

\medskip

{

\small

\textbf{Figure 3.}~Student's solution for Problem \ref{P-4}.

}

\end{center}

\medskip

\noindent{}We emphasize that mentioning the trivial statement (iii)$\Leftrightarrow$(iv)$\wedge$(v) but  not discussing the majority of the non-trivial statements illustrates how much freedom this type of exercises leaves to the students, although the setup looks like a `theorem construction kit'. The solution in Figure 3 is by far not satisfactory. Nevertheless it shows that the student analyzed the properties and recognized, in addition to the results from the lecture, that (iii)$\Leftrightarrow$(iv)$\wedge$(v) is trivially true. Again, we point out that recognizing what is trivially true is a core competence of a mathematician which is hard to pick up without discovery-oriented tasks.

\bigskip

\subsection{Another example}

In addition to the two courses mentioned above we used the following variants of our method also in other courses on real analysis. Notice, that in the course where Problem \ref{P-6} was given, compactness was defined via open covers and (i)$\Leftrightarrow$(iii) was thus a theorem.

\begin{adjustwidth}{1cm}{1cm}
\begin{h}\label{P-5} Let $f\colon\mathbb{R}\rightarrow\mathbb{R}$ be a map. Place the properties
\begin{compactitem}\vspace{3pt}
\item[(i)] $f$ is continuous,
 \vspace{2pt}
\item[(ii)] $f$ is differentiable,\vspace{3pt}
\item[(iii)] $f$ is continuously differentiable,\vspace{4pt}
\end{compactitem}
into the scheme below. Give a counterexample for each of the crossed out implications.\vspace{2pt}
\begin{center}\small
\begin{picture}(298,70)(30,0)
\put (0,0){\line (1,0){355}}
\put (0,65){\line (1,0){355}}
\put (0,0){\line (0,1){65}}
\put (355,0){\line (0,1){65}}

\put (103,35){\LARGE$\Longrightarrow$}
\put (110,20){\LARGE$\not\!\!\!\Longleftarrow$}

\put (228,35){\LARGE$\Longrightarrow$}
\put (235,20){\LARGE$\not\!\!\!\Longleftarrow$}

\end{picture}
\end{center}
\end{h}
\end{adjustwidth}

\medskip

\begin{adjustwidth}{1cm}{1cm}
\begin{h}\label{P-6}Let $(X,d)$ be a metric space and let $K\subseteq X$ be a subset. Place the properties
\begin{compactitem}\vspace{2pt}
\item[(i)] $K$ is compact,
 \vspace{2pt}
\item[(ii)] $K$ is closed and bounded,\vspace{2pt}
\item[(iii)] every sequence in $K$ has a convergent subsequence with limit in $K$,\vspace{4pt}
\end{compactitem}
into the scheme below. Give a counterexample for the crossed out implication.\vspace{2pt}
\begin{center}\small
\begin{picture}(298,70)(30,0)
\put (0,0){\line (1,0){355}}
\put (0,65){\line (1,0){355}}
\put (0,0){\line (0,1){65}}
\put (355,0){\line (0,1){65}}

\put (100,26){\LARGE$\Longleftrightarrow$}

\put (228,35){\LARGE$\Longrightarrow$}
\put (235,20){\LARGE$\not\!\!\!\Longleftarrow$}

\end{picture}
\end{center}
\end{h}
\end{adjustwidth}

\bigskip

The design in Problems \ref{P-5} and \ref{P-6} firstly again reduces the complexity as there now are only six respectively three conjectures. It also avoids `unwanted' answers like in Figure 2, where the student formed new conditions by connecting the given ones with logical operators. Our experiences with this design indicate that it can serve very well as an intermediate step before `the same' questions is given without pretyped scheme.


\bigskip

\section{Student's Feedback}\label{SEC:4}

\smallskip

During the course on Ordinary Differential Equations we gave in addition to Problem \ref{P-1} several other tasks of a discovery-oriented design. Among these were the request to construct a certain example (five times), and the question if an effect that was discovered earlier still can take place if the assumptions are changed (one time). The major part of the course's homework tasks ($12\times4=48$ questions) were rather classical textbook problems, that howewer included questions like `Prove [\,\dots\,]' or `Show [\,\dots\,]'.

\bigskip

At the end of the course we asked the participants to rate the discovery-oriented assignments by ticking one out of six boxes, e.g.,

\bigskip

\begin{center}
{not helpful}\;\;\;\;{\Large\Square\;\Square\;\Square\;\Square\;\Square\;\Square} \;\; {very helpful}
\end{center}

\bigskip

for the first question in Table 2. This design was used since in the German school system there are six grades (similar to A--F in many english-speaking countries). Below we quantified the answers with values 1 to 6 from the left to the right. The questions and the outcome read then as follows.

\bigskip

\begin{center}

{

\small

\textbf{Table 2.}\label{T-2}~Students' rating of our exercises.

}

\medskip

\begin{tabular}{cccc}
    \toprule
         & \hspace{-3pt}\phantom{$\displaystyle\sum$}N\hspace{-4pt}\phantom{$\displaystyle\sum$} & \hspace{-8pt}Mean\hspace{1pt} & SD\hspace{-3pt}\vspace{1pt}\\ \midrule\vspace{10pt}
\hspace{-5pt}\begin{minipage}{295pt}\small After completing the oral exam, please rate the discovery-\phantom{$\displaystyle\sum$}\\oriented problems (1=not helpful, 6=very helpful).\end{minipage}\hspace{5pt} & 11 & \hspace{-8pt}4.82&  1.17\hspace{-3pt} \\\vspace{10pt}
\hspace{-5pt}\begin{minipage}{295pt}\small Please rate how the discovery-oriented problems supported\\ your learning during the semester (1=not at all, 6=very much).\end{minipage}\hspace{5pt} & 11 &  \hspace{-8pt}4.36 &  1.36\hspace{-3pt} \\\vspace{10pt}
\hspace{-5pt}\begin{minipage}{295pt}\small Please compare the discovery-oriented problems to classical\\ones (1=like classical ones more, 6=prefer discovery-oriented).\end{minipage}\hspace{5pt} & 11 &  \hspace{-8pt}4.27 &  0.90\hspace{-3pt} \\\vspace{5pt}
\hspace{-5pt}\begin{minipage}{295pt}\small Would you like to work on discovery-oriented problems again\\in the future (1=don't want, 6=want very much)? \end{minipage}\hspace{5pt} & 11 &  \hspace{-8pt}4.91 &  0.94\hspace{-3pt} \\\bottomrule
  \end{tabular}

\end{center}

\bigskip

In their anonymous comments the students mentioned explicitly, that the discovery-oriented tasks encourage to study intensively (seven times mentioned), that they require to work independently (four times mentioned), and that they foster creativity (three times mentioned), e.g.:

\medskip

\begin{center}
\includegraphics[width=320pt]{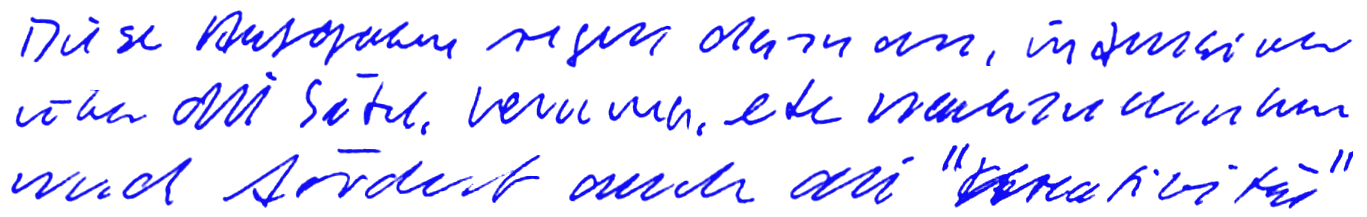}
\\\textquotedblleft{}These exercises encourage to think more intensively\\about the theorems, lemmas etc.~and foster the creativity.\textquotedblright{}
\end{center}

\medskip

\noindent{}The students mentioned that the tasks are difficult (four times mentioned), e.g.:

\begin{center}
\includegraphics[width=400pt]{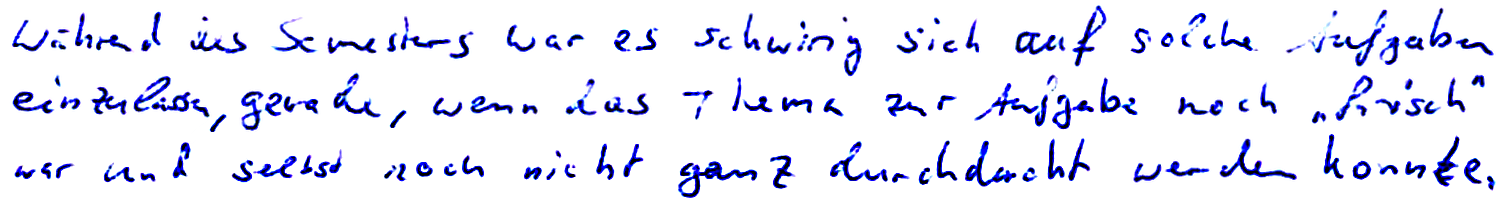}
\\\textquotedblleft{}During the semester it was difficult to engage with such problems,\\in particular when the topic was still very fresh and there had\\not been the time to think it through completely.\textquotedblright{}
\end{center}

\medskip

\noindent{}The answers however indicate that a suitable preparation by standard tasks helps to solve the discovery-oriented tasks successfully. Indeed, it was mentioned  explicitly that the tasks are not helpful without suitable preparation (three times mentioned), but helpful if prepared properly (four times mentioned), e.g.:

\medskip

\begin{center}
\includegraphics[width=400pt]{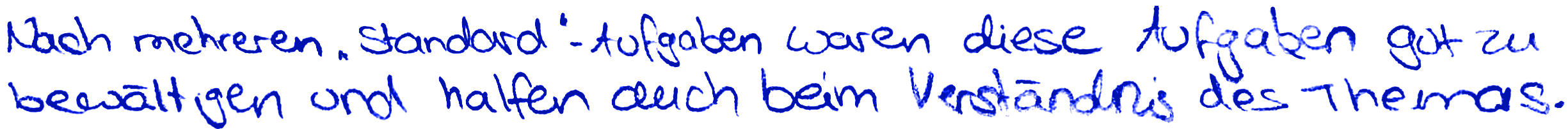}
\\\textquotedblleft{}After several standard exercises these exercises were\\well-manageable and helped also to grasp the topic.\textquotedblright{}
\end{center}

\medskip

\noindent{}In view of the formal assessment, the students mentioned that the discovery-oriented problems are helpful to prepare for the oral exam (three times mentioned) but it was remarked that with these problems it is difficult to obtain credit for grade improvement, e.g.:

\medskip

\begin{center}
\includegraphics[width=400pt]{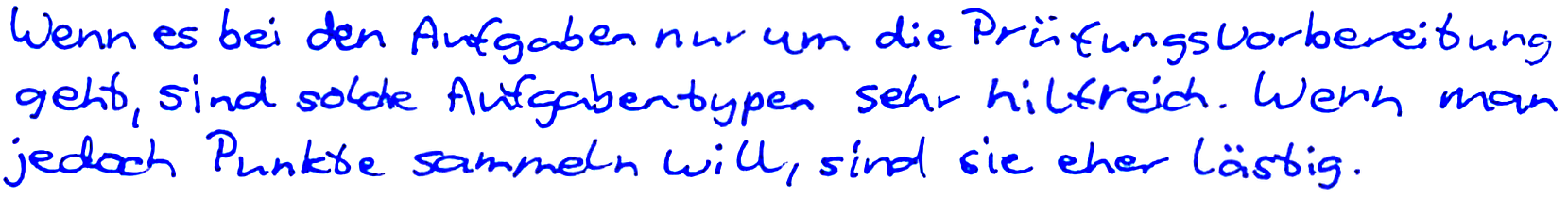}
\\\textquotedblleft{}If the problems are only for the preparation of the [oral] exam, then these types\\of tasks are very helpful. However, if one wants to collect credit, then they are rather annoying.\textquotedblright{}\label{COM}
\end{center}

\medskip

\noindent{}It was pointed out that both types of tasks are important as they address different levels of understanding. In particular, the students mentioned that the discovery-oriented tasks help to obtain a deeper understanding (seven times mentioned), e.g.: 

\medskip

\begin{center}
\includegraphics[width=400pt]{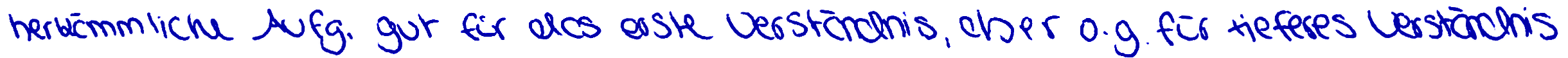}
\\\textquotedblleft{}Conventional exercises are good for the first understanding,\\but [\dots\,these tasks are good] for a deeper understanding.\textquotedblright{}
\end{center}


\bigskip
\section{Conclusion}\label{SEC:5}

\smallskip

Comparing Problem \ref{P-1} and Problem \ref{P-2} shows that Section \ref{SEC:3} gives rise to a general method which allows to transform certain classical exercises into discovery-oriented problems with a limited complexity. Our experiences and the feedback of our students suggest that the latter are suitable to serve as assessed homework, if (a) the  discovery-oriented problems are well-prepared by classical exercises and (b) additional study time and effort is credited directly or is motivated by highlighting the effect of the  problems on students' exam preparation.

\smallskip

We invite all teachers of mathematics to use and alterate the tasks that we explained in this article and to share their experiences with us.

\bigskip

\begin{center}

{\sc Acknowledgements} 

\end{center}

The authors would like to thank the anonymous referees for their careful work and valuable suggestions that improved this article significantly.


\normalsize\singlespace

\bibliographystyle{amsplain}

\end{document}